\documentclass{article}

\PassOptionsToPackage{numbers, compress}{natbib}

 \usepackage[main, final]{neurips_2025}

\usepackage[utf8]{inputenc} 
\usepackage[T1]{fontenc}    
\usepackage{hyperref}       
\usepackage{url}            
\usepackage{booktabs}       
\usepackage{amsfonts}       
\usepackage{nicefrac}       
\usepackage{microtype}      
\usepackage{xcolor}         

\usepackage{amssymb, amsmath}
\usepackage{stmaryrd}

\usepackage{comment}
\usepackage{bm}
\usepackage[textsize=tiny,
disable
]{todonotes}

\usepackage{tikz}
\usetikzlibrary{positioning}
\usepackage{pgfplots}
\pgfplotsset{compat=1.10}
\usepgfplotslibrary{groupplots,fillbetween}
\usetikzlibrary{spy, backgrounds}
\usetikzlibrary{shapes.geometric, arrows.meta}

\def\Dpartial#1#2{ \frac{\partial #1}{\partial #2} }

\def\DpartN#1#2#3{ \frac{\partial^#3 #1 }{ \partial #2^#3} }

\newcommand{\YC}[1]{#1}
\newcommand{\KC}[1]{#1}

\renewcommand{\Re}{\mathrm{Re}}
\newcommand{\En}{\mathrm{En}}
\newcommand{\De}{\mathrm{De}}

\newcommand{\MSE}{\mathrm{MSE}}
\newcommand{\dir}{\mathrm{dir}}

\newcommand{\cJ}{\mathcal{J}}

\newcommand{\cN}{\mathcal{N}}

\newcommand{\rR}{\mathbb{R}}

\newcommand{\bD}{\mathbf{D}}

\newcommand{\bN}{\mathbf{N}}

\newcommand{\bR}{\mathbf{R}}

\newcommand{\bQ}{\mathbf{Q}}

\newcommand{\bZ}{\mathbf{Z}}

\newcommand{\bq}{\mathbf{q}}
\newcommand{\bx}{\mathbf{x}}

\newcommand{\bz}{\mathbf{z}}

\newcommand{\f}{\mathbf{f}}
\newcommand{\bPhi}{\bm{\Phi}}
\newcommand{\bXi}{\bm{\Xi}}

\def\En{\mathrm{En}}
\def\De{\mathrm{De}}

\pgfplotsset{select coords between index/.style 2 args={
    x filter/.code={
        \ifnum\coordindex<#1\fi
        \ifnum\coordindex>#2\fi
    }
}}

\title{Latent space element method}

\author{%
  Seung Whan Chung\\
  Center for Applied Scientific Computing\\
  Lawrence Livermore National Laboratory\\
  Livermore, CA 94550 \\
  \texttt{chung28@llnl.gov} \\
  \And
  Youngsoo Choi\\
  Center for Applied Scientific Computing\\
  Lawrence Livermore National Laboratory\\
  Livermore, CA 94550 \\
  \texttt{choi15@llnl.gov} \\
  \And
  Christopher Miller\\
  Material Science Division\\
  Lawrence Livermore National Laboratory\\
  Livermore, CA 94550 \\
  \texttt{miller279@llnl.gov} \\
  \And
  H.\ Keo Springer\\
  Material Science Division\\
  Lawrence Livermore National Laboratory\\
  Livermore, CA 94550 \\
  \texttt{springer12@llnl.gov} \\
  \And
  Kyle T.\ Sullivan\\
  Material Science Division\\
  Lawrence Livermore National Laboratory\\
  Livermore, CA 94550 \\
  \texttt{sullivan34@llnl.gov} \\
}

\begin{document}

\maketitle

\begin{abstract}
\YC{How can we build surrogate solvers that train on small domains but scale to
larger ones without intrusive access to PDE operators? Inspired by the
Data-Driven Finite Element Method (DD-FEM) framework for modular data-driven
solvers, we propose the Latent Space Element Method (LSEM), an element-based
latent surrogate assembly approach in which a learned subdomain (“element”)
model can be tiled and coupled to form a larger computational domain. Each
element is a LaSDI latent ODE surrogate trained from snapshots on a local
patch, and neighboring elements are coupled through learned directional
interaction terms in latent space, avoiding Schwarz iterations and interface
residual evaluations. A smooth window-based blending reconstructs a global
field from overlapping element predictions, yielding a scalable assembled
latent dynamical system. Experiments on the 1D Burgers and Korteweg-de
Vries equations show that LSEM maintains predictive accuracy while scaling to
spatial domains larger than those seen in training. LSEM offers an
interpretable and extensible route toward foundation-model surrogate solvers
built from reusable local models.}

\end{abstract}

\section{Introduction}\label{sec:intro}

The rapid evolution of large language models
(LLMs)---including ChatGPT, Gemini, Claude, and others---has
demonstrated the transformative potential of foundation models
built from massive datasets and capable of general-purpose inference
across an extremely broad range of tasks.
The success of LLMs has naturally inspired interest
in developing foundational models for computational science,
prompting the question: Can scientific simulations benefit from
similarly universal, reusable, and highly scalable models?
However, foundational models in computational science differ
in several fundamental ways from their NLP counterparts~\cite{choi2025defining}.
Whereas LLMs prioritize semantic generalization over unstructured text,
scientific foundation models must reproduce governing physical laws,
maintain stability, and achieve controlled accuracy---all
while delivering orders-of-magnitude speedup over traditional numerical solvers.
This shift in objective reframes the design criteria: scientific foundational models
must be not only expressive but also physics-consistent, geometry-flexible,
and computationally efficient.
\par
In recent years, the scientific machine learning (SciML) community
has produced a wide range of promising data-driven and operator-learning approaches.
SINDy~\cite{brunton2016discovering} introduced sparse discovery of governing equations from data.
Neural operators such as the Fourier Neural Operator (FNO)~\cite{li2020fourier}
and the generalized neural operator framework~\cite{kovachki2023neural}
demonstrated strong accuracy and mesh-independent learning of PDE solution operators.
DeepONet~\cite{lu2019deeponet} also established theoretical universality for learning nonlinear operators.
While these approaches have led to impressive accelerations in PDE surrogate modeling,
most models are trained on fixed domains and thus lack the geometric flexibility
required for scalable scientific-foundation-model behavior.
\par
A parallel line of research focuses on data-driven finite element methods (DD-FEM)~\cite{choi2025defining},
which provide a modular and interpretable alternative.
DD-FEM leverages reduced-order models (ROMs) as reusable element-level solvers,
assembled via domain decomposition to enable prediction on arbitrary geometries or extended domains.
Applications include \YC{thermal transfer}~\cite{huynh2013static},
ROM-based lattice design~\cite{mcbane2021component, mcbane2022stress},
scalable time-dependent PDE solvers~\cite{chung2024train, chung2024scaled},
nonlinear manifold ROMs~\cite{diaz2024fast, zanardi2024scalable},
and improved Schwarz-based coupling strategies~\cite{wentland2025role}.
These works demonstrate that ROM-enabled finite elements form a promising blueprint
for geometry-flexible foundational models in computational science.
\par
However, all existing DD-FEM approaches remain intrusive,
requiring deep access to the underlying high-fidelity solver to extract snapshots,
residuals, and \YC{full order} operators.
This intrusiveness provides theoretical rigor but imposes significant overhead:
access to proprietary solvers may be restricted,
and constructing intrusive ROMs often requires substantial engineering
or even development of full-order models from scratch.
These limitations motivate a non-intrusive element framework
that preserves the geometric flexibility of DD-FEM
while removing the dependency on intrusive solvers.
\par
Non-intrusive models, however, introduce their own challenges.
Most SciML models are trained on fixed geometries and
therefore do not naturally generalize to unseen spatial configurations or larger domain sizes.
Some recent works attempt to address this through iterative coupling:
for example, Wang et al.~\cite{wang2022mosaic} used deep neural operators
combined with a Schwarz alternating method,
but the number of iterations grows rapidly with domain size.
Ouyang et al.~\cite{ouyang2025neural} embedded neural-operator elements (NOE)
into a classical FEM solver but did not demonstrate fully NOE-driven global predictions.
Thus, a practical and scalable non-intrusive DD-FEM pipeline
capable of arbitrary geometry assembly remains missing.
\par
In this work, we propose Latent Space Element Method (LSEM): a new DD-FEM
framework that uses a non-intrusive latent-space dynamics model---Latent Space
Dynamics Identification (LaSDI)~\cite{fries2022lasdi, tran2024weak,
bonneville2024gplasdi, park2024tlasdi, bonneville2024comprehensive,
anderson2025sequential, stephany2025rollout, stephany2025higher}---as an element solver, enabling
geometry- and scale-generalizable prediction without requiring iterative
Schwarz procedures or access to intrusive solvers.
\par
To achieve geometric and scale flexibility,
we deliberately design the framework around several key principles:
\begin{enumerate}
\item \textbf{Reusable latent-space element models}.
Each element is represented by a LaSDI model that learns a compact ODE
governing its local latent dynamics.
These element models are trained on reference subdomains
and reused across arbitrary global configurations.
\item \textbf{Smooth interface assembly}.
Neighboring elements are coupled through window-based interpolation
that enforces smooth transitions and continuity across overlapping regions,
allowing elements to be arranged in new geometries without retraining.
\item \textbf{Interaction dynamics learning}.
Latent variables of adjacent elements communicate
through a learned interaction model,
enabling information flow that mimics inter-element flux or coupling in classical PDE solvers.
\item \textbf{Scalable global evolution}.
The global system is advanced by assembling all local latent dynamics
into a single block-structured latent ODE system,
supporting extension to domains many times larger than those used in training.
\end{enumerate}
\par
Together, these design mechanisms allow latent space elements
to function analogously to classical finite elements:
trained once, assembled arbitrarily, and capable of accurate prediction
on geometries and scales never seen during training.
\par
The rest of the paper proceeds as follows.
Section~\ref{sec:method} introduces the methodology,
including latent dynamics identification, interface handling, and interaction modeling.
Section~\ref{sec:results} demonstrates the approach on Burgers' and KdV equations,
showcasing strong geometric scalability.
Section~\ref{sec:conclusion} concludes with future directions for foundational modeling
in computational science.

\section{Methodology}\label{sec:method}

\subsection{Latent Space Dynamics Identification}\label{subsec:lasdi}

Latent Space Dynamics Identification (LaSDI) provides the foundation
of our non-intrusive element model.
The goal of LaSDI is to learn a low-dimensional dynamical system
whose evolution accurately represents that of a high-dimensional physics state variable.
\KC{
Let $\bq(t)\in\rR^{N}$ denote the state variable defined on
fixed grid points in a single physics domain $\Omega$,
\begin{equation}
\bq(t) =
\begin{pmatrix}
q(\bx_1, t) & q(\bx_2, t) & \cdots & q(\bx_N, t)
\end{pmatrix},
\end{equation}
for a scalar variable $q(\bx, t)$.
The extension to vector variables should be straightforward.
Let the training data consist of a time series of snapshots,
}
\begin{equation}
\bQ =
\begin{pmatrix}
\bq(t_1) & \bq(t_2) & \cdots & \bq(t_{N_t})
\end{pmatrix}.
\end{equation}
\par
LaSDI begins by learning a nonlinear low-dimensional embedding
of the physics state using an autoencoder (AE).
The AE consists of an encoder $\En(\bq; \theta)$ and a decoder $\De(\bz; \theta)$,
parameterized by neural network weights $\theta$.
The encoder maps the high-dimensional state to a latent variable $\bz(t)\in \rR^{N_z}$
with $N_z\ll N$:
\begin{subequations}
\begin{equation}
\bz = \En(\bq; \theta),
\end{equation}
and the decoder reconstructs the physics state,
\begin{equation}
\bq = \De(\bz; \theta).
\end{equation}
\end{subequations}
Although LaSDI is agnostic to the AE architecture,
in this work we use fully connected neural networks
for both encoder and decoder to maintain simplicity and reproducibility.
\par
Once the physics trajectory is encoded,
LaSDI models the temporal evolution of the latent variable $\bz(t)$
as a first-order ordinary differential equation (ODE):
\begin{equation}\label{eq:ld-single}
\dot{\bz} = \bXi\bPhi(\bz).
\end{equation}
Here, $\bPhi\in\rR^{M}$ is a feature library consisting of
candidate nonlinear functions of $\bz$ (e.g., polynomials, interaction terms),
and $\bXi\in\rR^{N_z\times M}$ is the coefficient matrix to be learned.
This formulation parallels the structure of
Sparse Identification of Nonlinear Dynamics (SINDy)~\cite{brunton2016discovering},
but the dynamics identification is performed directly in latent space,
where the governing behavior of the reduced system tends
to be significantly simpler.
\par
\begin{figure*}[tbph]
    \input{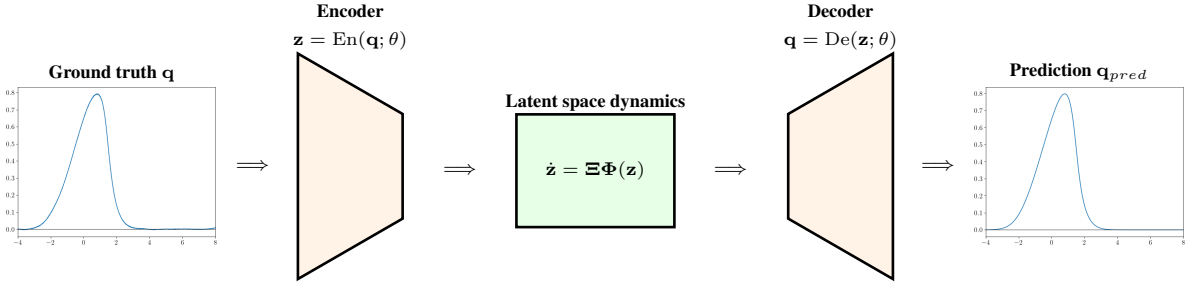}
    \caption{Schematic diagram of latent dynamics identification framework:
    The autoencoder maps the high-dimensional state $\bq$
    to a low-dimensional latent variable $\bz$,
    whose evolution is governed by the learned latent ODE $\dot{\bz} = \bXi\bPhi(\bz)$.
    The predicted latent trajectory is then decoded to
    reconstruct the physical state,
    enabling efficient reduced-order modeling of complex dynamical systems.}
    \label{fig:schematic}
\end{figure*}
A high-level schematic of the LaSDI process is illustrated in Figure~\ref{fig:schematic}.
In summary, the AE provides a nonlinear coordinate transformation
that compresses the physics state,
while the latent ODE captures its essential temporal evolution.
\par
The training process of LaSDI
consists of three primary components:
\begin{enumerate}
\item \textbf{Autoencoder (AE) loss}.\\
The encoded training trajectory
\begin{equation}
\bZ \equiv
\begin{pmatrix}
\bz(t_1) & \bz(t_2) & \cdots & \bz(t_{N_t})
\end{pmatrix}
= \En(\bQ)\in\rR^{N_z\times N_t}.
\end{equation}
is decoded back into physical space.
The AE reconstruction error is measured via the mean squared error (MSE):
\begin{equation}\label{eq:lasdi-Jae}
\cJ_{AE} = \MSE(\bQ, \De(\bZ)).
\end{equation}
\item \textbf{Latent dynamics (LD) loss}.\\
To ensure that the latent trajectory follows the ODE \eqref{eq:ld-single},
we enforce consistency between the learned right-hand side and
the finite-difference time derivative:
\begin{equation}\label{eq:lasdi-Jld}
\cJ_{LD} = \MSE(\bD_t\bZ^{\top}, \bPhi(\bZ)^{\top}\bXi^{\top}),
\end{equation}
where $\bD_t\in\rR^{N_t\times N_t}$ is a finite-difference operator for the time derivative.
\item \textbf{Regularization.}\\
Optionally, to avoid overfitting or ill-conditioned dynamics,
a regularization penalty may be imposed on $\bXi$:
\begin{equation}
\cJ_{reg} = \Vert \bXi \Vert^2.
\end{equation}
Though, in this study a different regularization loss is employed,
which will be introduced in subsequent sections.
\end{enumerate}
The total loss is a weighted combination:
\begin{equation}
\cJ = \cJ_{AE} + \alpha_{LD}\cJ_{LD} + \alpha_{reg}\cJ_{reg}.
\end{equation}
with hyperparameters $\alpha_{LD}$ and $\alpha_{reg}$
controlling the relative importance of dynamics fidelity and stability.
\par
For a single domain, the procedure above yields an efficient reduced model
capable of reconstructing and advancing the physics state.
To use LaSDI as an element model within our domain-decomposition framework,
however, additional components are required:
(i) methods for encoding and reconstructing subdomain states in overlapping regions,
and (ii) a learned mechanism for modeling inter-element interactions.
These extensions are introduced subsequently.

\subsection{Interface handling via smooth window interpolation}\label{subsec:window}

\begin{figure*}[tbph]
    \begin{tikzpicture}[
    font=\scriptsize,
    ]

\newcommand{\topleftsquare}[5]{%
  \draw[rounded corners=5pt, fill=#5, fill opacity=0.3, draw=none]
    (#1-#3/2, #2-#3/2) -- (#1+#3/2, #2-#3/2) [sharp corners] -- (#1+#3/2, #2+#3/2) -- (#1-#3/2, #2+#3/2) -- cycle;
  \ifx&#4&\else
    \node[draw=black, fill=white, line width=1.0, #4, scale=1.5] at (#1, #2) {};
  \fi
}

\newcommand{\topsquare}[5]{%
  \draw[rounded corners=5pt, fill=#5, fill opacity=0.3, draw=none]
    (#1-#3/2, #2+#3/2) -- (#1-#3/2, #2-#3/2) -- (#1+#3/2, #2-#3/2) [sharp corners] -- (#1+#3/2, #2+#3/2) -- cycle;
  \ifx&#4&\else
    \node[draw=black, fill=white, line width=1.0, #4, scale=1.5] at (#1, #2) {};
  \fi
}

\newcommand{\toprightsquare}[5]{%
  \draw[rounded corners=5pt, fill=#5, fill opacity=0.3, draw=none]
    (#1-#3/2, #2+#3/2) -- (#1-#3/2, #2-#3/2) [sharp corners] -- (#1+#3/2, #2-#3/2) -- (#1+#3/2, #2+#3/2) -- cycle;
  \ifx&#4&\else
    \node[draw=black, fill=white, line width=1.0, #4, scale=1.5] at (#1, #2) {};
  \fi
}

\newcommand{\leftsquare}[5]{%
  \draw[rounded corners=5pt, fill=#5, fill opacity=0.3, draw=none]
    (#1-#3/2, #2-#3/2) -- (#1+#3/2, #2-#3/2) -- (#1+#3/2, #2+#3/2) [sharp corners] -- (#1-#3/2, #2+#3/2) -- cycle;
  \ifx&#4&\else
    \node[draw=black, fill=white, line width=1.0, #4, scale=1.5] at (#1, #2) {};
  \fi
}

\newcommand{\centersquare}[5]{%
  \draw[rounded corners=5pt, fill=#5, fill opacity=0.3, draw=none]
    (#1-#3/2, #2-#3/2) -- (#1+#3/2, #2-#3/2) -- (#1+#3/2, #2+#3/2) -- (#1-#3/2, #2+#3/2) -- cycle;
  \ifx&#4&\else
    \node[draw=black, fill=white, line width=1.0, #4, scale=1.5] at (#1, #2) {};
  \fi
}

\newcommand{\rightsquare}[5]{%
  \draw[rounded corners=5pt, fill=#5, fill opacity=0.3, draw=none]
    (#1+#3/2, #2+#3/2) -- (#1-#3/2, #2+#3/2) -- (#1-#3/2, #2-#3/2) [sharp corners] -- (#1+#3/2, #2-#3/2) -- cycle;
  \ifx&#4&\else
    \node[draw=black, fill=white, line width=1.0, #4, scale=1.5] at (#1, #2) {};
  \fi
}

\newcommand{\bottomleftsquare}[5]{%
  \draw[rounded corners=5pt, fill=#5, fill opacity=0.3, draw=none]
    (#1+#3/2, #2-#3/2) -- (#1+#3/2, #2+#3/2) -- (#1-#3/2, #2+#3/2) [sharp corners] -- (#1-#3/2, #2-#3/2) -- cycle;
  \ifx&#4&\else
    \node[draw=black, fill=white, line width=1.0, #4, scale=1.5] at (#1, #2) {};
  \fi
}

\newcommand{\bottomsquare}[5]{%
  \draw[rounded corners=5pt, fill=#5, fill opacity=0.3, draw=none]
    (#1+#3/2, #2-#3/2) -- (#1+#3/2, #2+#3/2) -- (#1-#3/2, #2+#3/2) [sharp corners] -- (#1-#3/2, #2-#3/2) -- cycle;
  \ifx&#4&\else
    \node[draw=black, fill=white, line width=1.0, #4, scale=1.5] at (#1, #2) {};
  \fi
}

\newcommand{\bottomrightsquare}[5]{%
  \draw[rounded corners=5pt, fill=#5, fill opacity=0.3, draw=none]
    (#1+#3/2, #2+#3/2) -- (#1-#3/2, #2+#3/2) [sharp corners] -- (#1-#3/2, #2-#3/2) -- (#1+#3/2, #2-#3/2) -- cycle;
  \ifx&#4&\else
    \node[draw=black, fill=white, line width=1.0, #4, scale=1.5] at (#1, #2) {};
  \fi
}

\def\compsize{1.6}
\def\compgap{1.2}

\newcommand{\refsquare}[5]{%
  \draw[fill=#5, fill opacity=0.3,]
    (#1+#3/2, #2+#3/2) -- (#1-#3/2, #2+#3/2) -- (#1-#3/2, #2-#3/2) -- (#1+#3/2, #2-#3/2) -- cycle;
  \ifx&#4&\else
    \node[draw=black, fill=white, line width=1.0, #4, scale=1.5] at (#1, #2) {};
  \fi
}

\newcommand{\refAE}[5]{
    \refsquare{#1}{#2}{#3}{#4}{#5}

    \node at (#1 + \compsize/2 + 0.5, #2) {$\Longrightarrow$};

    \node[
        trapezium,
        trapezium angle=60,
        minimum width=\compsize cm,
        rotate=-90,
        fill=#5,
        fill opacity=0.3,
        draw=black,
        line width=1pt,
    ] at (#1 + \compsize + 0.5, #2) {};
    \node[
        trapezium,
        trapezium angle=60,
        minimum width=\compsize cm,
        rotate=90,
        fill=#5,
        fill opacity=0.3,
        draw=black,
        line width=1pt,
    ] at (#1 + \compsize + 1.5, #2) {};
}

\topleftsquare{0.5 * \compgap}{-0.5 * \compgap}{\compsize}{}{orange}
\topsquare{1.5 * \compgap}{-0.5 * \compgap}{\compsize}{star}{blue}
\topsquare{2.5 * \compgap}{-0.5 * \compgap}{\compsize}{}{orange}
\toprightsquare{3.5 * \compgap}{-0.5 * \compgap}{\compsize}{circle}{brown}

\leftsquare{0.5 * \compgap}{-1.5 * \compgap}{\compsize}{}{orange}
\centersquare{1.5 * \compgap}{-1.5 * \compgap}{\compsize}{diamond}{green}
\centersquare{2.5 * \compgap}{-1.5 * \compgap}{\compsize}{}{orange}
\rightsquare{3.5 * \compgap}{-1.5 * \compgap}{\compsize}{}{orange}

\leftsquare{0.5 * \compgap}{-2.5 * \compgap}{\compsize}{star}{blue}
\centersquare{1.5 * \compgap}{-2.5 * \compgap}{\compsize}{}{orange}
\centersquare{2.5 * \compgap}{-2.5 * \compgap}{\compsize}{circle}{brown}
\rightsquare{3.5 * \compgap}{-2.5 * \compgap}{\compsize}{}{orange}

\bottomleftsquare{0.5 * \compgap}{-3.5 * \compgap}{\compsize}{}{orange}
\bottomsquare{1.5 * \compgap}{-3.5 * \compgap}{\compsize}{}{orange}
\bottomsquare{2.5 * \compgap}{-3.5 * \compgap}{\compsize}{star}{blue}
\bottomrightsquare{3.5 * \compgap}{-3.5 * \compgap}{\compsize}{diamond}{green}

\node[font=\normalsize] at (2 * \compgap, -4.5 * \compgap) {Global configuration};

\def\refgap{1.2 * \compsize}
\refAE{7 * \compgap}{0 * \refgap}{\compsize}{}{orange}
\refAE{7 * \compgap}{-1 * \refgap}{\compsize}{star}{blue}
\refAE{7 * \compgap}{-2 * \refgap}{\compsize}{circle}{brown}
\refAE{7 * \compgap}{-3 * \refgap}{\compsize}{diamond}{green}
\node[anchor=south, font=\normalsize] at (7 * \compgap, 0.5*\compsize) {Subdomain type};
\node[anchor=south, font=\normalsize] at (7 * \compgap + \compsize + 1, 0.5*\compsize) {Autoencoder};

\pgfmathsetmacro{\arraywidth}{5.7}
\pgfmathsetmacro{\plotgap}{6.0} 

\begin{groupplot}[
    group style={
        group name=my plots,
        group size=1 by 1,
        xlabels at=edge bottom,
        horizontal sep=2cm,
        vertical sep=2.2cm,
    },
]    
\pgfplotsset{set layers=standard}

    \nextgroupplot[
        height=0.2*\arraywidth cm,
        width=\arraywidth cm,
        anchor=south west,
        ylabel={$W(\bx)$},
        xshift=-0.2cm, yshift=0.5cm,
        tick scale binop={\times},
        xtick=\empty,
        ytick=\empty,
        xticklabels={},
        yticklabels={},
        tickwidth=0pt,
        scale only axis,
        enlarge x limits={false, abs value = 5mm},
        ymin=0.,
    ]

        \addplot+ [
            line width=1.0,
            solid,
            mark=none,
            orange,
        ]
        table [
            x index=0, y index=4,
        ]{./data/windows.txt};

        \addplot+ [
            line width=1.0,
            solid,
            mark=none,
            blue,
        ]
        table [
            x index=1, y index=5,
        ]{./data/windows.txt};

        \addplot+ [
            line width=1.0,
            solid,
            mark=none,
            orange,
        ]
        table [
            x index=2, y index=6,
        ]{./data/windows.txt};

        \addplot+ [
            line width=1.0,
            solid,
            mark=none,
            brown,
        ]
        table [
            x index=3, y index=7,
        ]{./data/windows.txt};

\end{groupplot}

\draw[red, dashed, line width=1.5pt]
    (0.5*\compgap - \compsize/2 - 0.3, -0.5*\compgap) -- 
    (3.5*\compgap + \compsize/2 + 0.3, -0.5*\compgap);

\draw[->, red, line width=1.0pt, >=stealth]
    (0.5*\compgap - \compsize/2 - 0.3, -0.45*\compgap) 
    to[out=135, in=-135] 
    (0.45*\compgap - \compsize/2, 0.2*\compsize);
    
\end{tikzpicture}
    \caption{Schematic diagram of interface handling in the proposed LSEM framework:
    Local autoencoders are assigned to reference subdomain types,
    and global configurations are constructed by tiling these subdomains
    with overlapping regions. Smooth window functions $W(\bx)$
    blend the reconstructed subdomain solutions,
    ensuring continuity and smooth transitions across element interfaces.
    }
    \label{fig:interface-scheme}
\end{figure*}
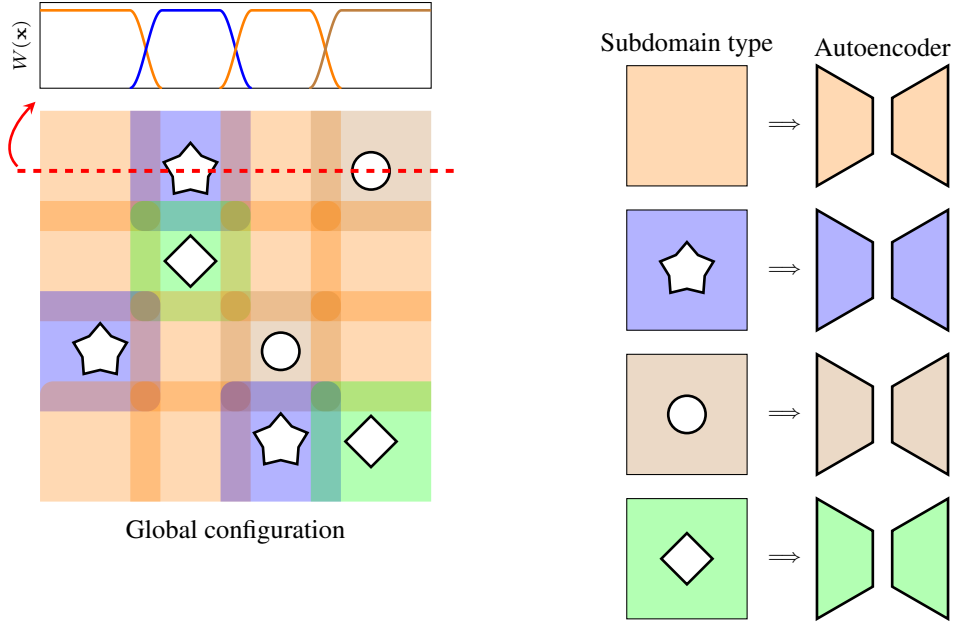
To extend LaSDI from single-domain modeling to
a geometry-flexible DD-FEM framework,
we must enable element-wise latent models to assemble into
a coherent global prediction.
Each element model is defined on a unit reference subdomain,
and global configurations are constructed by
tiling these reference subdomains across an arbitrary geometry
(see Figure~\ref{fig:interface-scheme}).
Neighboring subdomains overlap with each other,
the global state is reconstructed
via smooth window interpolation mechanism
that preserves continuity and smoothness across interfaces.
\par
Each reference subdomain type $r$ is assigned its own autoencoder,
enabling the framework to accommodate heterogeneous geometries or material regions.
For a subdomain $m$ of type $r(m)$,
the local state variable is encoded and decoded as
\begin{subequations}
\begin{equation}
\bz_{r} = \En_{r}(\bq_{r}; \theta_r)
\end{equation}
\begin{equation}
\bq_{r} = \De_{r}(\bz_{r}; \theta_r).
\end{equation}
\end{subequations}
Here, $\bq_m\in\rR^{N_m}$ denotes the local physics state associated with subdomain $m$.
In practice, these subdomains constitute the building blocks
from which arbitrarily large or geometrically complex configurations are assembled.
\par
Given a global state $\bq\in\rR^{N}$ defined over the assembled domain $\Omega$,
the local subdomain state $\bq_m$ is obtained through a restriction operator:
\begin{equation}\label{eq:Rq}
\bq_m = \bR_m\bq,
\end{equation}
where $\bR_m\in\rR^{N_m\times N}$ extracts the degrees of freedom
on the subdomain grid. Although written as a matrix for notational convenience,
this operation is implemented via direct slicing of $\bq$ based on subdomain indexing.
\par
After encoding, each latent variable $\bz_m$ evolves under
global latent dynamics model (introduced in Section~\ref{subsec:lsem-ld}).
The resulting reconstructed element states $\bq_m$
must then be merged into a single global field.
\par
Subdomain overlap is essential to enforce continuity and smoothness
across element boundaries without additional equation solving.
To merge the reconstructed subdomain states,
we form the global reconstruction by weighted superposition
of element contributions:
\begin{equation}\label{eq:Rtq}
\bq = \sum_m \bR_m^{\top}\mathrm{diag}(W_m(\bx_m))\bq_m.
\end{equation}
Here, $\bx_m\in\rR^{N_m\times d}$ denotes the coordinates of subdomain $m$, and
$W_m:\rR^d\to[0, 1]$ is a smooth window function that decays to zero at the
subdomain boundary.
\KC{
The window functions satisfy the partition of unity in the domain, i.e.
\begin{equation}
\sum_m W_m(\bx) = 1
\qquad
\forall \bx\in \Omega.
\end{equation}
}
This smooth decay ensures that overlapping subdomains
interpolate seamlessly, producing a globally continuous and differentiable
field without introducing artificial sharp transitions or oscillations.
\par
For example, on one-dimensional subdomain $[x_m, x_{m+1}]$ with overlap length $L_b$,
\begin{equation}
W_m(x) =
\begin{cases}
\frac{1}{2} - \frac{1}{2}\cos\left(\pi\frac{x - x_m}{L_b}\right) & x \in [x_m, x_m+L_b] \\
1 & x \in [x_m+L_b, x_{m+1}-L_b]\\
\frac{1}{2} - \frac{1}{2}\cos\left(\pi\frac{x - x_{m+1}}{L_b}\right) & x \in [x_{m+1}-L_b, x_{m+1}] \\
0 & \text{otherwise}.
\end{cases}
\end{equation}
This cosine-based window grants $C^1$-continuity
across element boundaries and aligns with established practices
in overlapping domain decomposition and partition-of-unity finite element methods.

\subsection{Modeling of interaction dynamics between elements}\label{subsec:lsem-ld}

To assemble local LaSDI element models into a coherent global system,
it is necessary to model not only the internal dynamics of each element
but also the interactions between neighboring elements.
Although each subdomain is equipped with its own encoder, decoder,
and latent evolution law, the resulting elements do not evolve in isolation;
their latent variables must exchange information in a manner analogous
to flux or coupling terms in traditional domain decomposition
or finite volume methods.
The purpose of this section is to describe
how we introduce these interactions directly in latent space.
\par
For simplicity, we assume that all subdomains share the same
latent dynamics formulation, regardless of subdomain type.
This assumption reflects the physical expectation
that the governing PDE is the same everywhere and that differences
in local geometry can be captured by the element-specific autoencoders
rather than by modifying the latent ODE structure.
The assumption is not essential, however,
and can be relaxed if distinct physics or material models
require heterogeneous dynamics across elements.
\par
\begin{figure*}[tbph]
    \begin{tikzpicture}[
    font=\scriptsize,
    ]

\newcommand{\topleftsquare}[5]{%
  \draw[rounded corners=5pt, fill=#5, fill opacity=0.3, draw=none]
    (#1-#3/2, #2-#3/2) -- (#1+#3/2, #2-#3/2) [sharp corners] -- (#1+#3/2, #2+#3/2) -- (#1-#3/2, #2+#3/2) -- cycle;
  \ifx&#4&\else
    \node[draw=black, fill=white, line width=1.0, #4, scale=1.5] at (#1, #2) {};
  \fi
}

\newcommand{\topsquare}[5]{%
  \draw[rounded corners=5pt, fill=#5, fill opacity=0.3, draw=none]
    (#1-#3/2, #2+#3/2) -- (#1-#3/2, #2-#3/2) -- (#1+#3/2, #2-#3/2) [sharp corners] -- (#1+#3/2, #2+#3/2) -- cycle;
  \ifx&#4&\else
    \node[draw=black, fill=white, line width=1.0, #4, scale=1.5] at (#1, #2) {};
  \fi
}

\newcommand{\toprightsquare}[5]{%
  \draw[rounded corners=5pt, fill=#5, fill opacity=0.3, draw=none]
    (#1-#3/2, #2+#3/2) -- (#1-#3/2, #2-#3/2) [sharp corners] -- (#1+#3/2, #2-#3/2) -- (#1+#3/2, #2+#3/2) -- cycle;
  \ifx&#4&\else
    \node[draw=black, fill=white, line width=1.0, #4, scale=1.5] at (#1, #2) {};
  \fi
}

\newcommand{\leftsquare}[5]{%
  \draw[rounded corners=5pt, fill=#5, fill opacity=0.3, draw=none]
    (#1-#3/2, #2-#3/2) -- (#1+#3/2, #2-#3/2) -- (#1+#3/2, #2+#3/2) [sharp corners] -- (#1-#3/2, #2+#3/2) -- cycle;
  \ifx&#4&\else
    \node[draw=black, fill=white, line width=1.0, #4, scale=1.5] at (#1, #2) {};
  \fi
}

\newcommand{\centersquare}[5]{%
  \draw[rounded corners=5pt, fill=#5, fill opacity=0.3, draw=none]
    (#1-#3/2, #2-#3/2) -- (#1+#3/2, #2-#3/2) -- (#1+#3/2, #2+#3/2) -- (#1-#3/2, #2+#3/2) -- cycle;
  \ifx&#4&\else
    \node[draw=black, fill=white, line width=1.0, #4, scale=1.5] at (#1, #2) {};
  \fi
}

\newcommand{\rightsquare}[5]{%
  \draw[rounded corners=5pt, fill=#5, fill opacity=0.3, draw=none]
    (#1+#3/2, #2+#3/2) -- (#1-#3/2, #2+#3/2) -- (#1-#3/2, #2-#3/2) [sharp corners] -- (#1+#3/2, #2-#3/2) -- cycle;
  \ifx&#4&\else
    \node[draw=black, fill=white, line width=1.0, #4, scale=1.5] at (#1, #2) {};
  \fi
}

\newcommand{\bottomleftsquare}[5]{%
  \draw[rounded corners=5pt, fill=#5, fill opacity=0.3, draw=none]
    (#1+#3/2, #2-#3/2) -- (#1+#3/2, #2+#3/2) -- (#1-#3/2, #2+#3/2) [sharp corners] -- (#1-#3/2, #2-#3/2) -- cycle;
  \ifx&#4&\else
    \node[draw=black, fill=white, line width=1.0, #4, scale=1.5] at (#1, #2) {};
  \fi
}

\newcommand{\bottomsquare}[5]{%
  \draw[rounded corners=5pt, fill=#5, fill opacity=0.3, draw=none]
    (#1+#3/2, #2-#3/2) -- (#1+#3/2, #2+#3/2) -- (#1-#3/2, #2+#3/2) [sharp corners] -- (#1-#3/2, #2-#3/2) -- cycle;
  \ifx&#4&\else
    \node[draw=black, fill=white, line width=1.0, #4, scale=1.5] at (#1, #2) {};
  \fi
}

\newcommand{\bottomrightsquare}[5]{%
  \draw[rounded corners=5pt, fill=#5, fill opacity=0.3, draw=none]
    (#1+#3/2, #2+#3/2) -- (#1-#3/2, #2+#3/2) [sharp corners] -- (#1-#3/2, #2-#3/2) -- (#1+#3/2, #2-#3/2) -- cycle;
  \ifx&#4&\else
    \node[draw=black, fill=white, line width=1.0, #4, scale=1.5] at (#1, #2) {};
  \fi
}

\def\compsize{1.6}
\def\compgap{3}

\topleftsquare{0.5 * \compgap}{-0.5 * \compgap}{\compsize}{}{orange}
\topsquare{1.5 * \compgap}{-0.5 * \compgap}{\compsize}{star}{blue}
\topsquare{2.5 * \compgap}{-0.5 * \compgap}{\compsize}{}{orange}

\leftsquare{0.5 * \compgap}{-1.5 * \compgap}{\compsize}{}{orange}
\centersquare{1.5 * \compgap}{-1.5 * \compgap}{\compsize}{diamond}{green}
\centersquare{2.5 * \compgap}{-1.5 * \compgap}{\compsize}{}{orange}

\leftsquare{0.5 * \compgap}{-2.5 * \compgap}{\compsize}{star}{blue}
\centersquare{1.5 * \compgap}{-2.5 * \compgap}{\compsize}{}{orange}
\centersquare{2.5 * \compgap}{-2.5 * \compgap}{\compsize}{circle}{brown}


\pgfmathsetmacro{\cx}{1.5 * \compgap}
\pgfmathsetmacro{\cy}{-1.5 * \compgap}

\draw[->, >=stealth, line width=1pt] 
    (0.5*\compgap + \compsize/2, -0.5*\compgap - \compsize/2) -- 
    (\cx - \compsize/2, \cy + \compsize/2)
    node[midway, below left] {$\bPhi(\bz_{11})\bXi_{tl}$};

\draw[->, >=stealth, line width=1pt] 
    (1.5*\compgap, -0.5*\compgap - \compsize/2) -- 
    (\cx, \cy + \compsize/2)
    node[midway, right, yshift=5pt,] {$\bPhi(\bz_{12})\bXi_{t}$};

\draw[->, >=stealth, line width=1pt] 
    (2.5*\compgap - \compsize/2, -0.5*\compgap - \compsize/2) -- 
    (\cx + \compsize/2, \cy + \compsize/2)
    node[midway, right] {$\bPhi(\bz_{13})\bXi_{tr}$};

\draw[->, >=stealth, line width=1pt] 
    (0.5*\compgap + \compsize/2, -1.5*\compgap) -- 
    (\cx - \compsize/2, \cy)
    node[midway, above] {$\bPhi(\bz_{21})\bXi_{l}$};

\draw[->, >=stealth, line width=1pt] 
    (2.5*\compgap - \compsize/2, -1.5*\compgap) -- 
    (\cx + \compsize/2, \cy)
    node[midway, above] {$\bPhi(\bz_{23})\bXi_{r}$};

\draw[->, >=stealth, line width=1pt] 
    (0.5*\compgap + \compsize/2, -2.5*\compgap + \compsize/2) -- 
    (\cx - \compsize/2, \cy - \compsize/2)
    node[midway, left] {$\bPhi(\bz_{31})\bXi_{bl}$};

\draw[->, >=stealth, line width=1pt] 
    (1.5*\compgap, -2.5*\compgap + \compsize/2) -- 
    (\cx, \cy - \compsize/2)
    node[midway, right, yshift=-5pt,] {$\bPhi(\bz_{32})\bXi_{b}$};

\draw[->, >=stealth, line width=1pt] 
    (2.5*\compgap - \compsize/2, -2.5*\compgap + \compsize/2) -- 
    (\cx + \compsize/2, \cy - \compsize/2)
    node[midway, right] {$\bPhi(\bz_{33})\bXi_{br}$};

\node[anchor=south] at (\cx, \cy - \compsize/2) {$\bPhi(\bz_{22})\bXi_{D}$};

\end{tikzpicture}
    \caption{Schematic diagram of interaction dynamics modeling between neighboring latent-space elements:
    Each element’s latent variable evolves according to its internal dynamics
    together with contributions from neighboring elements,
    represented by directional coefficient blocks (e.g., $\bXi_t$, $\bXi_l$, $\bXi_{br}$).
    This structure enables the model to capture information flow
    across elements and assemble a coherent global latent dynamics operator.
    }
    \label{fig:interaction-scheme}
\end{figure*}
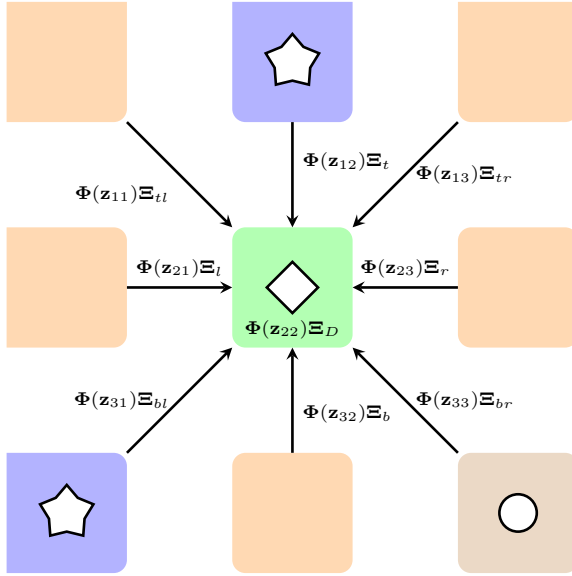
Consider the element located at position $(2,2)$ in the configuration 
illustrated in Figure~\ref{fig:interaction-scheme}.
Its latent variable $\bz_{22}$ evolves according
to a combination of internal dynamics and interaction contributions
from its neighboring elements:
\begin{equation}
\dot{\bz}_{22}
=
\underbrace{\bXi_D \bPhi(\bz_{22})}_{\text{internal dynamics}}
+
\underbrace{\f[\bz_{22}; \bN(\bz_{22})]}_{\text{interaction dynamics}},
\end{equation}
where $\bXi_D$ is the coefficient matrix governing the internal latent dynamics,
$\bPhi(\cdot)$ is the latent feature library, and
\begin{equation}
\bN(\bz_{22}) = \{\bz_{11}, \bz_{12}, \bz_{13}, \bz_{21}, \bz_{23}, \bz_{31}, \bz_{32}, \bz_{33} \}
\end{equation}
denotes the set of latent variables associated with the eight neighboring elements.
The function $\f$ characterizes how neighboring states influence the evolution of $\bz_{22}$.
\par
One of the simplest and most effective ways to construct $\f$
is to model the contribution from each neighbor using another LaSDI-type term,
producing
\begin{equation}\label{eq:upwind}
\f = \sum_{\bz_n\in \bN(\bz_{22})} \bXi_{\dir(\bz_n; \bz_{22})}\bPhi(\bz_n),
\end{equation}
where the direction index $\dir(\bz_{n};\bz_{22}) \in \{tl, t, tr, l, r, bl, b, br\}$
specifies whether the neighboring element lies top-left, top, top-right, and so on relative to element $(2, 2)$.
This directional indexing allows the learned interactions to distinguish between,
for example, left-to-right advection and top-to-bottom propagation.
In hyperbolic systems with a dominant flow direction,
the corresponding directional matrix (e.g., $\bXi_r$ or $\bXi_l$)
naturally acquires a larger magnitude,
enabling the latent model to mimic upwind flux behavior.
\par
Of course, more complex systems may require richer interaction formulations,
particularly where flow direction is not dominant or varies significantly over time.
In such cases, the interaction may also include the contribution
of the center element itself:
\begin{equation}\label{eq:lf-flux}
\f = \sum_{\bz_n\in \bN(\bz_{22})} \bXi_{\dir(\bz_n; \bz_{22})}\bPhi(\bz_n) + \bXi_{\dir(\bz_{22}; \bz_n)}\bPhi(\bz_{22}).
\end{equation}
This bidirectional formulation enhances
expressive power and enables the latent dynamics model
to emulate numerical flux mechanisms such as the Lax-Friedrichs flux
used in finite volume methods.
\par
The matrices $\bXi_D$ and all directional matrices $\bXi_{dir}$ together
serve as building blocks for the global dynamics operator.
When the elements are assembled into a larger configuration,
these blocks are arranged into a global coefficient matrix $\bXi$
with the up-wind model \eqref{eq:upwind},
\begin{equation}
\bXi=
\begin{pmatrix}
\bXi_D & \bXi_r & & & \bXi_b & \bXi_{br} \\
\bXi_l & \bXi_D & \bXi_r & & \bXi_{bl} & \bXi_{b} & \bXi_{br} \\
& \bXi_l & \bXi_D & \bXi_r & & \bXi_{bl} & \bXi_{b} & \bXi_{br} & \cdots\\
& & \bXi_l & \bXi_D & & & \bXi_{bl} & \bXi_{b} \\
\bXi_t & \bXi_{tl} & & & \bXi_D & \bXi_r \\
\bXi_{tr} & \bXi_{t} & \bXi_{tl} & & \bXi_l & \bXi_D & \bXi_r \\
& \bXi_{tr} & \bXi_{t} & \bXi_{tl} & & \bXi_l & \bXi_D & \bXi_r & \cdots\\
& & \bXi_{tr} & \bXi_{t} & & & \bXi_{l} & \bXi_{D} \\
& & \vdots & & & & \vdots & & \ddots \\
\end{pmatrix}
\in \rR^{16N_z\times 16M}.
\end{equation}
This block-structured organization encodes all information
about local dynamics, neighbor interactions, and element connectivity.
\par
For the latent dynamics, all subdomain latent variables are merged into
a global latent variable,
\begin{equation}
\bz = (\bz_{11}, \bz_{12}, \bz_{13}, \bz_{14}, \bz_{21}, \ldots, \bz_{44}) \in \rR ^{16N_z}.
\end{equation}
\KC{
The global feature library $\bPhi(\bz)$ takes
element-wise feature libraries,
\begin{equation}
\bPhi(\bz)=
(\bPhi(\bz_{11}), \bPhi(\bz_{12}), \ldots, \bPhi(\bz_{44}))
\in \rR^{16M}.
\end{equation}
}
Combined with the global feature library $\bPhi(\bz)$,
the resulting global dynamics model takes the form
\begin{equation}\label{eq:global-ld}
\dot{\bz} = \bXi\bPhi(\bz).
\end{equation}
\par
In this way, latent interactions are learned directly from data,
enabling information transfer across element boundaries
without requiring intrusive flux evaluation or iterative Schwarz-type procedures.
This latent-space coupling is what ultimately allows the element models
to operate collectively as a large-scale, geometry-flexible surrogate
consistent with the underlying physical dynamics.

\subsection{Reparameterization technique}\label{subsec:reparam}

When LaSDI is applied as an element model within a large, assembled DD-FEM system,
the latent variables must evolve robustly under time integration.
However, when overfitted, standard autoencoders often exhibit sensitivity
in the latent space: a small perturbation in $\bz$ can decode
into a disproportionately large perturbation in the physical state $\bq$.
Due to this sensitivity, even small integration errors can excite unphysical modes,
causing inaccurate reconstructions.
To mitigate these issues, we incorporate a reparameterization technique
inspired by a component of variational autoencoders.
\par
The key idea is to expose the decoder
to perturbed latent variables during training,
forcing the autoencoder to learn representations
that remain stable under small deviations.
Instead of evaluating the reconstruction loss
using the exact encoded latent trajectory $\bZ$,
we add stochastic perturbations and train the decoder
to accurately reconstruct the physics state from these noisy inputs:
\begin{equation}\label{eq:reparam}
\cJ_{AE} = \MSE(\bQ, \De(\bZ + \delta\bZ)).
\end{equation}
The noise $\delta\bZ=(\delta\bz(t_1), \delta\bz(t_2), \ldots, \delta\bz(t_{N_t}))\in\rR^{N_z\times N_t}$
consists of instantaneous noise terms drawn from a normal distribution,
\begin{equation}
\delta\bz(t) \sim \cN(0, \epsilon(t)^2).
\end{equation}
Rather than prescribing a fixed noise magnitude,
we adaptively determine the standard deviation $\epsilon(t)$
based on the dynamics residual of the latent ODE model.
Using the encoded training trajectory $\bZ$,
we evaluate the mismatch between the finite-difference time derivative
and the latent ODE prediction:
\begin{equation}
\mathrm{RMS}(\bD_t\bZ^{\top} - \bPhi(\bZ)^{\top}\bXi^{\top}).
\end{equation}
This residual provides a natural measure of how much local error
one may expect during time integration.
We therefore define the time-dependent noise level as
\begin{equation}
\epsilon(t) = \beta\frac{t}{\Delta t}\mathrm{RMS}(\bD_t\bZ^{\top} - \bPhi(\bZ)^{\top}\bXi^{\top}),
\end{equation}
where $\beta$ is a tunable hyperparameter controlling the aggressiveness of regularization.
\par
This construction produces two desirable effects.
First, because the perturbation grows linearly in time,
it mirrors the amplification that naturally occurs in numerical integration
of stable linear systems,
ensuring that the decoder is trained to handle the types of deviations
that will occur during prediction.
Second, the reparameterization effectively smooths the latent manifold
by preventing the decoder from overfitting to a thin submanifold of
exact encoded values.
As a result, the latent representation becomes more robust to integration drift
and modeling imperfections.
\par
Although the technique resembles the noise injection used in
variational autoencoders, we emphasize that we do not include
a KL-divergence regularization.
\KC{
While KL divergence regularize the latent space distribution
by adding an additional constraint,
numerical experiments showed that the noise injection was sufficient
enough to smooth the latent space without over-constraining the optimization.
}
Instead, the noise improves stability by broadening the latent distribution
sampled during training,
enforcing that nearby latent states decode to physically consistent fields.
In practice, this significantly reduces the risk of unphysical artifacts.

\subsection{Regularization in training}\label{subsec:regularization}

A critical requirement for reliable latent-space dynamics is
that the learned ODE remains numerically stable when integrated forward in time,
particularly in large assembled configurations
where modeling errors may propagate across many interacting elements.
While the reparameterization technique in Section 2.4 improves robustness
by regularizing the autoencoder,
the stability of the latent dynamics themselves must also be enforced
during training.
To achieve this, we incorporate regularization strategies
that penalize unstable behavior in the learned coefficient matrices.
\par
When the feature library $\bPhi$ contains only linear terms,
the global latent dynamics reduce to a linear system
whose stability can be assessed directly through eigenvalue analysis.
In this case, we compute the eigenvalues $\Lambda(\bXi)$
of the operator in \eqref{eq:global-ld} and penalize any with positive real part:
\begin{equation}\label{eq:Jreg1}
\cJ_{reg} = \mathrm{ReLU}(\Re(\Lambda(\bXi)))^2.
\end{equation}
This term enforces that the learned operator does not introduce latent modes
that grow exponentially, ensuring that time integration remains stable
for both training and unseen configurations.
\par
For nonlinear feature libraries---which are often necessary for capturing
complex PDE behavior---direct eigenvalue analysis is no longer applicable,
since the stability of the full nonlinear system depends on $\bz(t)$
and therefore varies in time. In this setting, we adopt an energy-based
regularization approach that evaluates the instantaneous rate of change
of the latent-state energy along the encoded training trajectory.
For each snapshot $\bz(t_k)$, the energy growth rate of
the latent dynamics \eqref{eq:global-ld} is
\begin{equation}
\frac{d}{dt}\left(\frac{\Vert\bz(t_k)\Vert^2}{2}\right)
= \bz(t_k)^{\top}\dot{\bz}(t_k)
= \bz(t_k)^{\top}\bXi\bPhi(\bz(t_k)).
\end{equation}
Averaging these contributions across all snapshots yields the regularization term
\begin{equation}\label{eq:Jreg2}
\cJ_{reg} = \frac{1}{N_t}\sum_{k}^{N_t}\mathrm{ReLU}\left(\bz(t_k)^{\top}\bXi\bPhi(\bz(t_k))\right).
\end{equation}
Penalizing positive values of this quantity discourages
latent dynamics that increase energy along the training trajectories,
promoting stability without requiring explicit linearization.
\par
Both regularization forms---eigenvalue-based for linear systems and
energy-based for nonlinear systems---serve the same purpose:
to prevent the learned latent dynamics from drifting into unstable regimes.
By integrating these regularizers into the full training objective,
the LaSDI models learn not only to match the observed dynamics
but also to evolve in a manner consistent with the stability characteristics
of the underlying physical system.


\subsection{Overall training procedure}

Because the method is designed to generalize across arbitrary global geometries,
the training data consist of multiple configurations constructed
from the same set of reference elements.
Each configuration provides a time series of high-dimensional physics states,
from which we extract subdomain-level information and
assemble the corresponding global latent trajectories.
\par
Suppose we are given $N_{train}$ training simulations,
each defined on a domain composed of $R$ reference element types.
For simulation $n$, the global physics trajectory $\bQ^{(n)}$
is decomposed into subdomain states $\bq_m^{(n)}$
using the restriction rule in \eqref{eq:Rq}.
Each subdomain state is then encoded into its latent representation via
\begin{equation}
\bz_m^{(n)} = \En_{r(m)}(\bq_m^{(n)};\theta_r),
\end{equation}
producing a set of latent trajectories that together
form the global latent trajectory $\bZ^{(n)}$.
From these trajectories, we construct the element-wise feature evaluations $\bPhi(\bZ^{(n)})$
and assemble the global dynamics matrix $\bXi$ following the block structure
described in \eqref{eq:global-ld}.
This process ensures that each training configuration contributes
information to all internal and directional dynamics blocks,
enabling the model to learn consistent interactions
regardless of the final assembly geometry.
\par
With the latent trajectories and global dynamics operators in hand,
the latent dynamics loss is evaluated across all training simulations,
\begin{equation}
\cJ_{LD} = \frac{1}{N_{train}}\sum_n^{N_{train}}\MSE(\bD_t\bZ^{(n)\top}, \bPhi(\bZ^{(n)})^{\top}\bXi^{\top}),
\end{equation}
which enforces that the learned operator correctly predicts
latent time derivatives across heterogeneous configurations.
\par
In parallel, we evaluate the autoencoder reconstruction loss
using the reparameterized formulation introduced in Section~\ref{subsec:reparam}.
For each subdomain and each training simulation, we add noise $\delta\bZ_m^{(n)}$
to the latent trajectory and decode to obtain reconstructed states.
The corresponding loss is
\begin{equation}
\cJ_{AE} = \frac{1}{N_{train}}\sum_n^{N_{train}}\MSE(\bQ_m^{(n)}, \De_{r(m)}(\bZ_m^{(n)} + \delta \bZ_m^{(n)})),
\end{equation}
where the summation spans all subdomains $m$ of all training simulations.
Notably, we evaluate reconstruction quality locally at the subdomain level
rather than at the global level.
This choice ensures that each decoder learns
to capture fine-scale structure within its subdomain
while leaving the global assembly to be handled
by the windowing mechanism of Section~\ref{subsec:window}.
\par
To promote dynamic stability, we incorporate the regularization terms
described in Section~\ref{subsec:regularization}.
For each training simulation, the contribution is
\begin{equation}
\cJ_{reg} = \sum_{n}^{N_{train}}\mathrm{ReLU}(\Re(\Lambda(\bXi^{(n)})))^2,
\end{equation}
or its nonlinear alternative \eqref{eq:Jreg2}, depending on the feature library in use.
\par
Finally, all components are combined into the global training objective,
\begin{equation}
\cJ = \alpha_{AE}\cJ_{AE} + \alpha_{LD}\cJ_{LD} + \alpha_{reg}\cJ_{reg},
\end{equation}
where the weights $\alpha_{AE}$, $\alpha_{LD}$, $\alpha_{reg}$
balance reconstruction fidelity, dynamical accuracy, and stability.
The optimization variables consist of the autoencoder parameters $\theta_r$
for each reference element type $r=1,\ldots,R$
and the full collection of coefficient matrices
for internal and directional latent dynamics.
Training proceeds using gradient-based optimization,
during which the model learns
(i) how to compress and reconstruct subdomain states,
(ii) how to evolve latent variables consistently with the underlying PDE,
and (iii) how to couple latent dynamics across elements
in a geometry-agnostic manner.
\par
Through this unified training process,
the resulting LaSDI elements become reusable building blocks
capable of assembling into large-scale, arbitrarily shaped configurations
while preserving both stability and predictive fidelity.

\section{Results}\label{sec:results}

\todo[inline]{
Youngsoo:
We need to include wall-clock (speed-up vs FOM), training time (probably
rough order), and plot for inference time vs number of elements, etc.
I think linear model for latent dynamic identification is used for both
numerical experiments? Therefore, the eigen regularization is used? I do not
think these choices are mentioned anywhere in Results section. Let's mention
them.
How about adding one ablation, i.e., Overlap width $L_b$ sensitivity: set it to be small/medium/large and see how the performance of LSEM changes.
Show the numerical effect of the reparameterization, $\beta$?
}




\subsection{One-dimensional Burgers' equation}\label{subsec:burgers}

We first consider the one-dimensional inviscid Burgers' equation
as a minimal test case to demonstrate that the proposed LSEM framework
can accurately reproduce local dynamics
while generalizing to significantly larger domains than those in training data.
This example is intentionally simple in physics
but nontrivial in its use of minimal training data for
scaled-up prediction capability.
\par
The governing equation is
\begin{equation}
\Dpartial{q}{t} + q\Dpartial{q}{x} = 0,
\end{equation}
posed on a spatial domain with homogeneous Dirichlet boundary conditions.
All training simulations are performed on the domain $x\in[-4, 8]$
over the time interval $t\in[0, 1]$, discretized with $N_x=2048$
grid points and timestep $\Delta t=10^{-3}$.
The full-order model consists of a first-order upwind finite-difference
spatial discretization combined with first-order implicit backward Euler
time integration.
\par
\begin{figure*}[tbph]
    \includegraphics[width=\textwidth]{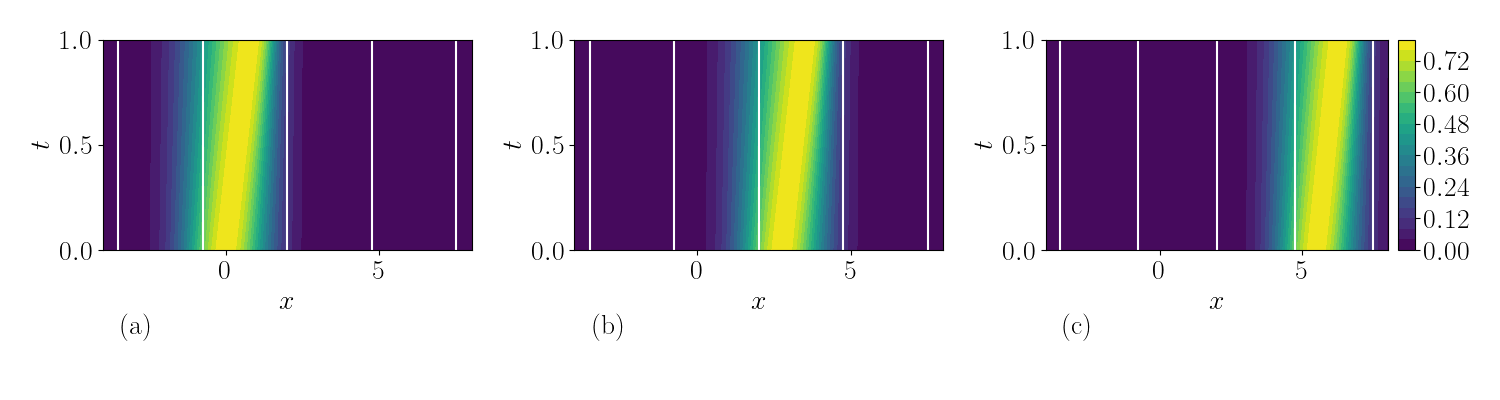}
    \caption{
        Four-element training datasets for one-dimensional Burgers' equation.
        White lines indicate element interfaces.
    }
    \label{fig:burgers-train}
\end{figure*}
To ensure that the latent interaction model captures directional propagation,
the training domain is partitioned into four overlapping elements.
Each element spans a region of length $L_e=\frac{14}{3}$,
and contains 639 grid points, including those in the overlap.
The overlap width is set to $L_b=1$.
The midpoints of overlapping regions, as shown in Figure~\ref{fig:burgers-train},
serve as effective element interfaces.
\begin{figure*}[tbph]
    \includegraphics[width=\textwidth]{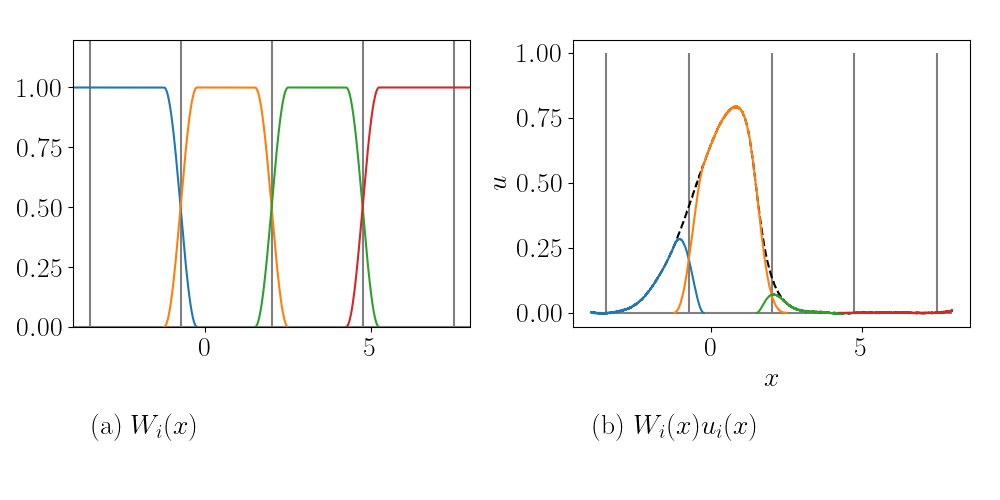}
    \caption{Solution composition on the four-element domain:
    (a) window functions $W_i(x)$ at each element $i=1,\ldots,4$; and
    (b) element-wise contribution $W_i(x)u_i(x)$,
    with black dashed line showing the global solution $u(x) = \sum_{i}W_i(x)u_i(x)$.
    Gray lines mark the element interfaces.
    }
    \label{fig:burgers-sol-composition}
\end{figure*}
Figure~\ref{fig:burgers-sol-composition} illustrates how the window functions $W_m(x)$
blend the subdomain reconstructions to form a smooth global solution.
\par
The initial condition for each training simulation is a shifted Gaussian pulse,
\begin{equation}
q(x, t=0) = A\exp\left(-\frac{(x-x_n)^2}{w^2}\right),
\end{equation}
with amplitude $A=0.8$ and width $w=1$.
The shift $x_n$ varies across training runs according to
\begin{equation}
x_n = 1.25 + nL_e, \qquad n=1,2,3.
\end{equation}
This design aligns the peak of the initial condition
at the same relative location within each element—specifically
at a distance $0.75$ from the left interface of its host element---while shifting
its global coordinate.
Consequently, the subdomain-level snapshots become identical
across all training datasets,
enabling efficient training of the autoencoder and element-wise latent dynamics.
On the other hand, the global latent trajectories differ due to the shifted positions,
ensuring that all blocks of the global dynamics matrix receive
sufficient training data and allowing the model to reproduce
the same initial pulse at an arbitrary cell location.
\par
The latent-space model for each element uses dimension $N_z=5$,
and all elements share a single autoencoder architecture:
an encoder with layer sizes $[639, 100, 30, 5]$
and a decoder with the reversed structure.
Softplus activations are used throughout.
Inter-element dynamics follow the one-way single-block formulation \eqref{eq:upwind},
which is sufficient for this simple solution
where information propagates predominantly from left to right.
\par
Training is performed using the SOAP optimizer for $2000$ epochs
with learning rate $3\times10^{-3}$,
and loss weights $\alpha_AE=\alpha_{LD}=1$, $\alpha_{reg}=0$.
No stability regularization is required due to the simplicity of the dynamics
and the linear-dominant feature set.
\par
\begin{figure*}[tbph]
    \includegraphics[width=\textwidth]{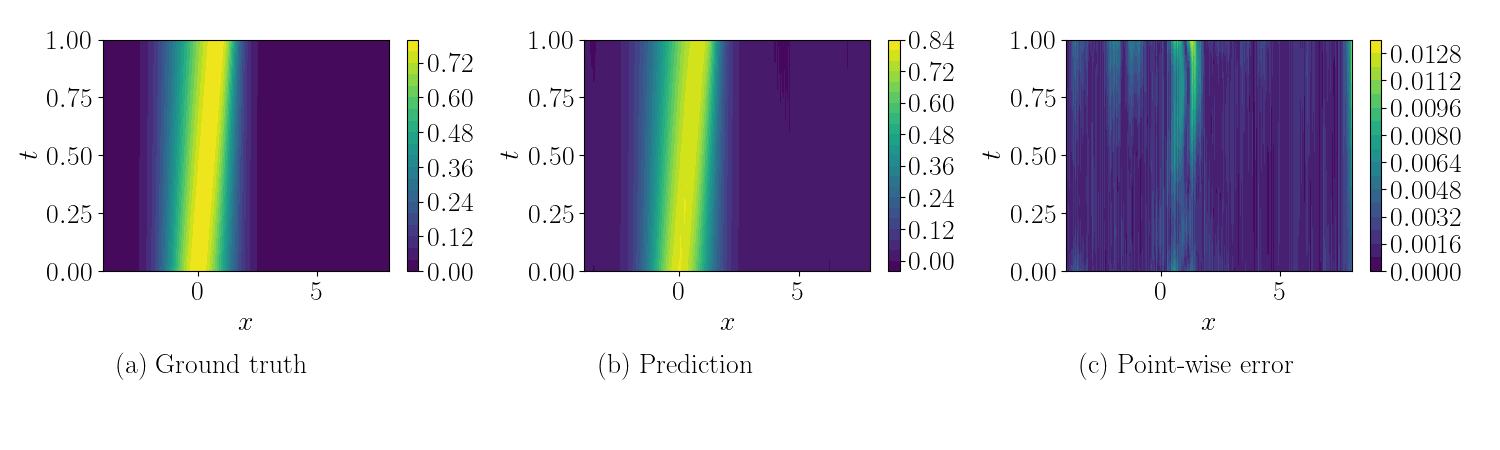}
    \caption{
    Reproductive prediction of Burgers' equation on the four-element training domain:
    (a) the full-order model (ground truth);
    (b) the LSEM prediction;
    (c) the pointwise absolute error over space-time; and
    (d) a final-time snapshot comparing the global prediction (solid line)
    and full-order solution (dashed line).
    The learned latent element model accurately reproduces
    the evolution of the initial pulse,
    with errors remaining small throughout the domain.
    }
    \label{fig:burgers-train-pred}
\end{figure*}
Figure~\ref{fig:burgers-train-pred} shows a reproductive prediction for the first training dataset.
The reconstructed global solution matches the full-order model closely,
with a relative $L_2$ error of $0.51\%$,
and the window-based reconstruction produces a smooth, continuous global field.
\par
To evaluate scalability, we test the trained model on a domain three times
larger than the training domain, $x\in[-4,32]$,
decomposed into 12 elements.
The initial condition consists of three Gaussian pulses:
\begin{equation}
q(x, t=0) = \sum_{k=1}^3 A\exp\left(-\frac{(x-x_k)^2}{w^2}\right),
\end{equation}
with centers $x_k=1.25 + X_kL_e$ and $X_k=1, 6, 9$.
This arrangement places pulses at widely separated positions
across the extended domain.
\begin{figure*}[tbph]
    \includegraphics[width=\textwidth]{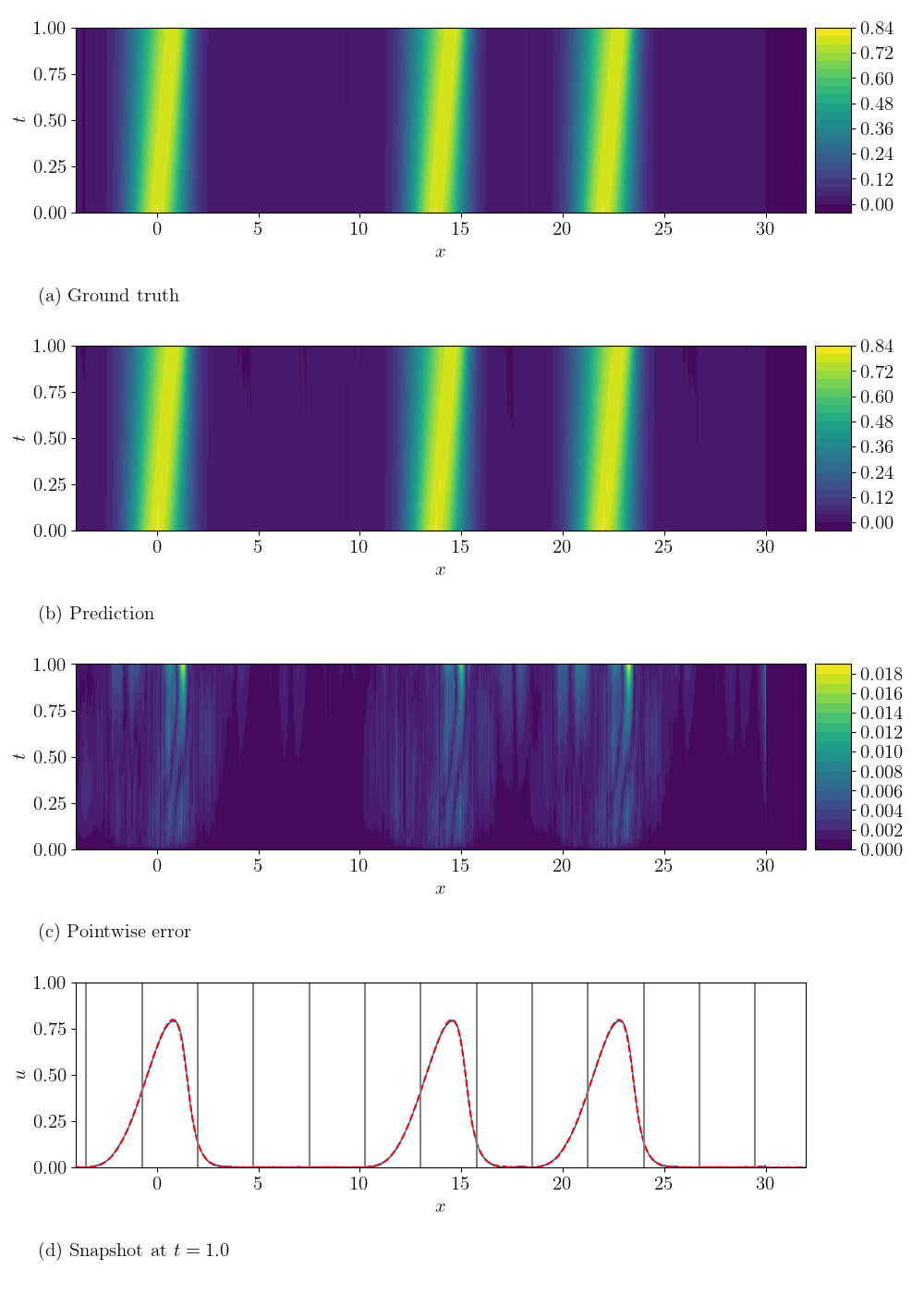}
    \caption{
    Prediction of Burgers' equation on a scaled-up 12-element domain:
    (a) ground truth;
    (b) LSEM prediction; and
    (c) pointwise absolute error over space-time; and
    (d) a final-time snapshot comparing the global prediction (solid line)
    and full-order solution (dashed line).
    The model generalizes from four training elements to twelve elements,
    maintaining accurate pulse shapes and positions with only minor localized errors.
    }
    \label{fig:burgers-scaledup-pred}
\end{figure*}
Figure~\ref{fig:burgers-scaledup-pred} displays the resulting prediction.
The framework successfully propagates all three pulses with high accuracy,
achieving a relative $L_2$ error of $0.41\%$.
\KC{Including the global interpolation \eqref{eq:Rtq},
this prediction was about $32$ times faster than the full-order model.}
No artificial transitions or oscillations were found in the prediction.

\subsection{One-dimensional Korteweg-De Vries equation}\label{subsec:kdv}

We next consider the one-dimensional Korteweg-De Vries (KdV) equation,
which presents a more challenging test due to its dispersive wave dynamics
and nonlinear wave-wave interaction.
Unlike Burgers' equation, where a single advecting pulse
provides relatively simple dynamics,
the KdV system exhibits soliton propagation and collision phenomena
that are highly dependent on initial conditions.
This example therefore serves to evaluate both the expressive capacity
of the latent dynamics model and the robustness of the LSEM assembly
to varied local wave behaviors.
\par
The governing equation is
\begin{equation}
\Dpartial{q}{t} + 6q\Dpartial{q}{x} + \DpartN{q}{x}{3}= 0,
\end{equation}
posed with periodic boundary conditions.
Training simulations are performed on the domain $x\in[-10, 30]$,
over the time interval $t\in[0, 1]$,
using $N_x=2000$ spatial grid points and timestep $\Delta t=10^{-3}$.
The full-order model uses a 1st-2nd order summation-by-parts (SBP)
finite-difference scheme for spatial derivatives
and a fourth-order explicit Runge-Kutta time integrator.
\par
To construct the training configurations,
the domain is divided into four overlapping elements,
each of length $L_e=10$ and containing $600$ grid points.
The overlap width is set to $L_b=2$,
ensuring smooth reconstruction across subdomain boundaries
via the window functions shown in Figure~\ref{fig:burgers-sol-composition}.
\par
The initial condition for each training simulation consists of
a random arrangement of solitons,
designed to produce a diverse set of interaction scenarios.
We begin with the two-soliton solution of the KdV equation~\cite{gardner1967method},
\begin{equation}
q(x, t=0) = \frac{6}{\cosh(x)^2},
\end{equation}
and construct a composite initial condition
by randomly assigning solitons to each cell:
\begin{subequations}\label{eq:kdv-ic}
\begin{equation}
q(x, t=0) = \sum_{k=1}^4A_k\frac{6}{\cosh(x-x_k)^2},
\end{equation}
where amplitude $A_k$
\begin{equation}
A_k\sim\mathrm{Bernoulli}(\frac{1}{2}),
\end{equation}
randomly turns each soliton ``on'' or ``off'',
and the center location $x_k$
\begin{equation}
x_k\sim \cN\!\!\left[\left(k - \frac{3}{2}\right)L_e, (0.1L_e)^2\right],
\end{equation}
is sampled near the center of each element with small spatial variance.
\end{subequations}
This process generates a wide variety of soliton interactions---merging,
splitting, and phase-shifting---which in turn provide rich training data
for both internal dynamics and inter-element coupling.
A total of 100 such training datasets are generated;
trivial cases with no solitons are removed.
\begin{figure*}[tbph]
    \includegraphics[width=\textwidth]{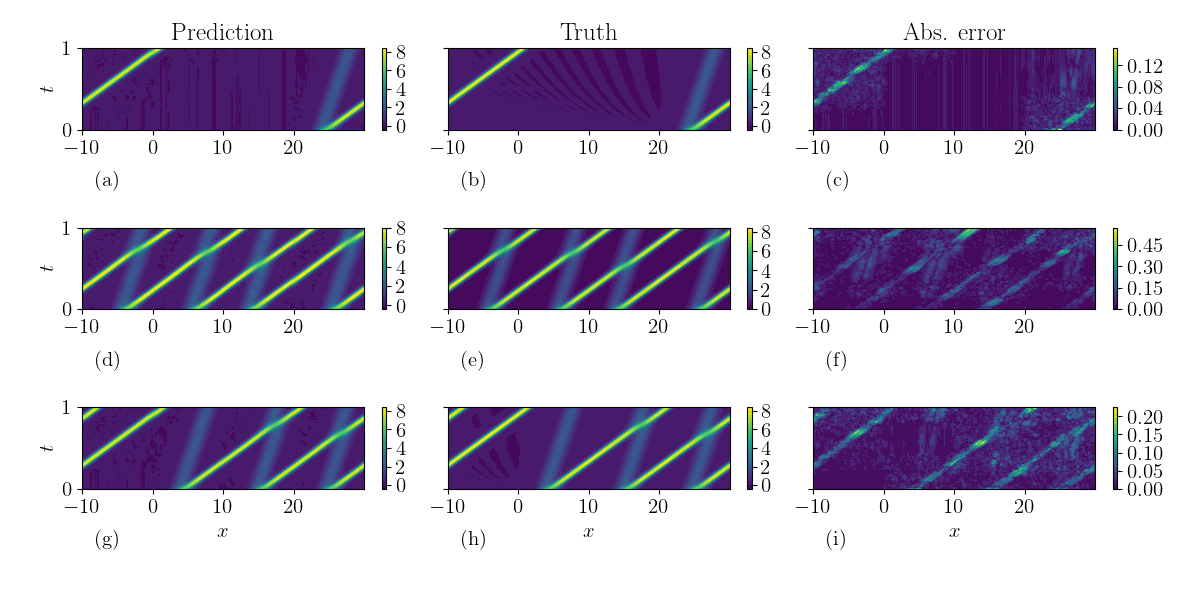}
    \caption{Prediction on the 42th, 46th, and 54th KdV training datasets.
    Columns show the model prediction, ground truth, and absolute error,
    demonstrating accurate reproduction of soliton propagation and interactions.
    }
    \label{fig:kdv-pred-train}
\end{figure*}
Representative examples are shown in Figure~\ref{fig:kdv-pred-train},
illustrating the diversity of collision patterns the model must capture.
\par
Each element uses a latent space of dimension $N_z=7$,
with a shared autoencoder architecture across all elements.
The encoder consists of layers $[600, 300, 100, 100, 30, 30, 10, 7]$,
and the decoder mirrors this structure.
$\tanh$ activations are used, as they provide smoother
latent representations for the dispersive waveforms characteristic of KdV.
Interaction dynamics again follow the one-way formulation \eqref{eq:upwind},
which we found sufficiently expressive for capturing soliton propagation
across elements.
\par
Training is performed using the SOAP optimizer for $10,000$ epochs
with learning rate $3\times10^{-3}$.
Loss weights are set to $\alpha_{AE}=\alpha_{reg}=1$ and $\alpha_{LD}=10$,
reflecting the increased importance of enforcing accurate latent dynamics
for systems with complex temporal behavior.
Regularization is particularly important here
to mitigate long-term phase drift and latent-mode instability.
\par
\begin{figure*}[tbph]
    \includegraphics[width=\textwidth]{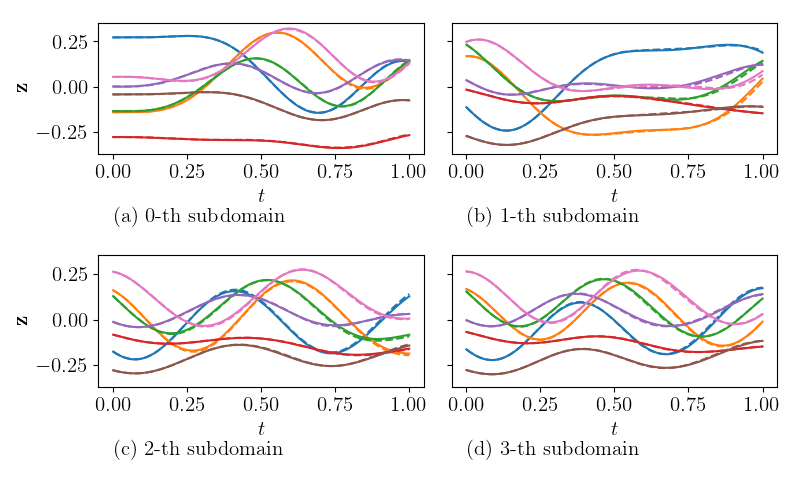}
    \caption{Latent-space trajectories for the 54th KdV training dataset:
    Predicted trajectories are compared with the encoded ground-truth
    trajectories (dashed lines), showing close agreement across all elements.
    }
    \label{fig:kdv-ld}
\end{figure*}
Across all training datasets, the autoencoder achieves a mean relative
$L_2$ reconstruction error of $0.76\%$,
indicating that the latent variables capture the soliton structures
with sufficient fidelity. The reproductive prediction
for a representative training case is shown in Figure~\ref{fig:kdv-pred-train},
where the model accurately captures the timing and position of soliton collisions.
The latent trajectories for the 54th dataset are shown in Figure~\ref{fig:kdv-ld},
illustrating that the predicted latent evolution matches
the encoded ground truth closely despite the highly nonlinear interactions.
The relative $L_2$ error for reproductive prediction across the dataset is $5.64\%$,
with most error arising from small phase mismatches in the peak soliton locations.
\par
\begin{figure*}[tbph]
    \includegraphics[width=\textwidth]{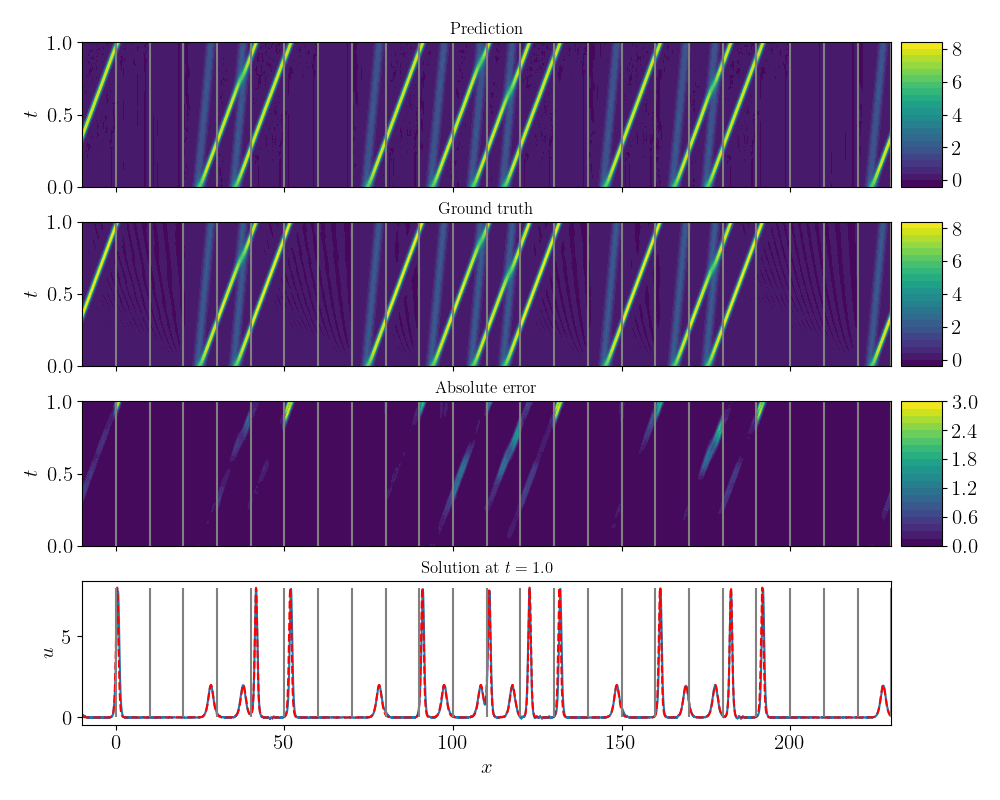}
    \caption{Prediction of KdV equation on a scaled-up 24-element domain.
    }
    \label{fig:kdv-scaledup}
\end{figure*}
To evaluate scalability, we test the trained model on
a domain six times larger than the training domain $x\in[-10, 230]$,
divided into 24 elements.
Initial conditions are generated using the same stochastic procedure \eqref{eq:kdv-ic}.
Figure~\ref{fig:kdv-scaledup} shows the prediction for
a representative configuration.
The model successfully propagates multiple solitons across long distances,
handling repeated interactions and collisions
despite never having been trained on domains of this size.
\KC{Including the global interpolation \eqref{eq:Rtq},
this prediction was about $5500$ times faster than the full-order model.}
The relative $L_2$ error of $11.7\%$ is dominated by slight phase drift
in the largest-amplitude soliton, a common difficulty in reduced models
for dispersive systems,
but the qualitative dynamics---including collision order and amplitude preservation---remain accurate.
\par
Overall, this experiment demonstrates that the proposed framework can
capture rich nonlinear wave interactions and maintain predictive capability
across domains far larger than those used for training,
highlighting its potential for scalable simulation of complex dynamics.

\section{Conclusion}\label{sec:conclusion}

This work introduced Latent Space Element Method---a non-intrusive
domain decomposition finite element framework
that employs Latent Space Dynamics Identification (LaSDI)
as a reusable element model for scalable prediction of dynamical systems.
By learning element-level latent dynamics directly from data and assembling
them through smooth window blending and latent-space interaction models,
the framework enables accurate prediction on geometries and domain sizes
far beyond those used during training.
\par
The numerical experiments highlight the strengths of this approach.
For Burgers' equation, the model achieved high accuracy despite extremely
limited training data, demonstrating the ability to ``train small, model big''.
For the more complex Korteweg-De Vries equation,
the framework successfully reproduced nonlinear soliton propagation
and collision patterns across domains many times larger than the training
configurations, with remaining discrepancies primarily due to mild phase drift.
These results underscore the method's capacity to combine non-intrusive learning,
geometric flexibility, and stable latent dynamics into a cohesive global predictor.
\par
While promising, the present study is limited to one-dimensional systems.
Extending the methodology to higher dimensions, more complex PDEs,
and problems with strong multi-scale or chaotic behavior will require
advances in encoder/decoder architectures, interaction modeling,
and stability control.
Future directions include two- and three-dimensional applications,
physics-informed feature libraries, adaptive element refinement,
and uncertainty quantification.
Overall, the proposed framework provides a scalable and flexible
foundation for developing next-generation,
geometry-generalizable models in computational science.

\begin{ack}
This work was performed under the auspices of the U.S. Department of Energy by
Lawrence Livermore National Laboratory under contract DE-AC52-07NA27344 and was
supported by Laboratory Directed Research and Development funding under project
24-SI-004. LLNL release number: LLNL-JRNL-2014444.
\end{ack}

\appendix

\bibliographystyle{elsarticle-num} 
\bibliography{reference.bib}

\end{document}